\documentclass[pre,preprint,amsmath,superscriptaddress,endfloats]{revtex4}

\usepackage{graphicx}

\begin{document}

\title{Synchronization in networks of networks: 
The onset of coherent collective behavior in systems of interacting
populations of heterogeneous oscillators}

\author{Ernest Barreto}
\email[Email address: ]{ebarreto@gmu.edu}
\affiliation{Department of Physics \& Astronomy, The Center for Neural Dynamics,
and The Krasnow Institute for Advanced Study, George Mason University, Fairfax
Virginia 22030, USA}

\author{Brian Hunt}
\affiliation{Department of Mathematics and Institute for Physical Science and Technology,
University of Maryland, College Park, Maryland 20742, USA}

\author{Edward Ott}
\affiliation{Institute for Research in Electronics and Applied Physics,
Department of Physics, and Department of Electrical and Computer Engineering,
University of Maryland, College Park, Maryland 20742, USA}

\author{Paul So}
\affiliation{Department of Physics \& Astronomy, The Center for Neural Dynamics,
and The Krasnow Institute for Advanced Study, George Mason University, Fairfax
Virginia 22030, USA}


\begin{abstract}
The onset of synchronization in networks of networks is investigated.
Specifically, we consider networks of interacting
phase oscillators in which the set of
oscillators is composed of several distinct populations. The oscillators
in a given population are heterogeneous in that their natural frequencies
are drawn from a given distribution, and each population has
its own such distribution. The coupling among the oscillators is global, however,
we permit the coupling strengths between the members of different populations to
be separately specified.
We determine the critical condition for the onset
of coherent collective behavior, and develop the illustrative case in
which the oscillator frequencies are drawn from a set of
(possibly different) Cauchy-Lorentz distributions.
One motivation is drawn from neurobiology, in which the collective dynamics of
several interacting populations of oscillators (such as excitatory and inhibitory
neurons and glia) are of interest.
\end{abstract}

\date{\today}

\maketitle

In recent years, there has been considerable interest in networks of interacting systems.
Researchers have found that an appropriate description of such systems involves
an understanding of both the dynamics of the individual oscillators and the connection
topology of the network. Investigators studying the latter have found that many complex
networks have a modular structure involving motifs \cite{motifs}, communities \cite{community,community2},
layers \cite{layers}, or clusters \cite{cluster}. For example, recent work has shown that as many kinds
of networks (including isotropic homogeneous networks and a class of scale-free
networks) transition to full synchronization,
they pass through epochs in which well-defined synchronized communities appear and interact
with one another \cite{community2}. Knowledge of this structure, and the dynamical behavior
it supports, informs our understanding of 
biological \cite{ZhouKurths}, social \cite{GirvanNewman}, and technological networks \cite{milo2}.

Here we consider the onset of coherent collective behavior in similarly structured systems
for which the dynamics of the individual oscillators is simple. In
seminal work, Kuramoto analyzed a mathematical model that
illuminates the mechanisms by which synchronization arises in a large set of
globally-coupled phase oscillators \cite{kur84}. An important feature of
Kuramoto's model is that the oscillators are heterogeneous in their frequencies.
And, although these mathematical results assume global coupling, they have been fruitfully
applied to further our understanding of systems of fireflies, arrays of Josephson
junctions, electrochemical oscillators, and many other cases \cite{apps}.

In this work, we study systems of several interacting Kuramoto systems, i.e., networks of
interacting {\it populations} of phase oscillators. Our motivation is drawn
not only from the applications listed above (e.g., an amusing application might be interacting populations of
fireflies, where each population inhabits a separate tree), but also from other biological
systems. Rhythms arising from coupled cell populations
are seen in many of the body's organs (including the heart, the pancreas,
and the kidneys, to name but a few), all of which interact. For example, the circadian rhythm
influences many of these systems. We draw additional motivation from
consideration of how neuronal tissue is organized. Although we do not consider
neuronal systems specifically in this paper, we note that heterogeneous ensembles
of neurons often exhibit a ``network-of-networks" topology. At the cellular level, populations of particular
kinds of neurons (e.g., excitatory neurons) interact not only among themselves, but
also with populations of other distinct neuronal types (e.g., inhibitory neurons).
At a higher level of organization, various anatomical regions of the brain interact
with one another as well \cite{ZhouKurths}.

Although our network is simple, it is novel in that we include heterogeneity
at several levels. Each population consists of a collection of phase oscillators
whose natural frequencies are drawn at random from a given distribution. To allow for
heterogenity at the population level, we let each population have its own such frequency
distribution. In addition, our system is heterogeneous at the coupling level as well: 
we consider global coupling such that the coupling strengths between the members of
different populations can be separately specified.

The assumpion of global (but population-weighted) coupling permits an analytic
determination of the critical condition for the onset of coherent collective behavior,
as we will show. While this assumption may not strictly apply in some applications, 
our results provide insight into the behavior
of networks of similar topology even if the connectivity is less than global.

We begin by specifying our model and deriving our main results. We then discuss
several illustrative examples.

Consider first a two-species Kuramoto model. We label the separate
species $\theta$ and $\phi$ and assume that there are $N_\theta$
and $N_\phi$ such oscillators in each population, respectively. Thus,
the system equations are
\begin{eqnarray}
\frac{d\theta_i}{dt} &=& \omega_{\theta i} + \eta \frac{k_{\theta \theta}}{N_\theta} \sum_{j=1}^{N_\theta} sin(\theta_{j} - \theta_i - \alpha_{\theta \theta})
+ \eta \frac{k_{\theta \phi}}{N_\phi} \sum_{j=1}^{N_\phi} sin(\phi_{j} - \theta_i - \alpha_{\theta \phi}), \nonumber \\
\frac{d\phi_i}{dt} &=& \omega_{\phi i} + \eta \frac{k_{\phi \theta}}{N_\theta} \sum_{j=1}^{N_\theta} sin(\theta_{j} - \phi_i - \alpha_{\phi \theta})
+ \eta \frac{k_{\phi \phi}}{N_\phi} \sum_{j=1}^{N_\phi} sin(\phi_{j} - \phi_i - \alpha_{\phi \phi}). \nonumber
\end{eqnarray}
Here, $\eta$ is an overall coupling strength, the $\alpha$'s provide additional heterogeneity in the coupling functions,
and the matrix
\begin{equation}
{\bf K}=\left(
\begin{array}{cc}
k_{\theta \theta} & k_{\theta \phi} \\
k_{\phi \theta} & k_{\phi \phi}
\end{array}
\right)
\label{Kmatrix}
\end{equation}
defines the connectivity among the populations \cite{Montbrio}.

For arbitrarily many different species, let $\sigma$ range over the various population
symbols $\theta, \phi, \dots$ with the understanding that
depending on the context, $\sigma$ may represent either a label (when subscripted)
or a variable. Thus, we have
\[
\frac{d\sigma_i}{dt} = \omega_{\sigma i}
+ \sum_{\sigma'} \left[
\eta \frac{k_{\sigma \sigma'}}{N_{\sigma'}} \sum_{j=1}^{N_{\sigma'}} sin(\sigma'_{j} - \sigma_i
- \alpha_{\sigma \sigma'})
\right].
\]
The $\omega_{\sigma i}$ are the constant natural frequencies of the oscillators when uncoupled,
and are distributed according to a set of time-independent distribution functions
$G_\sigma(\omega_\sigma)$ \cite{restreponote}.

We define the usual Kuramoto order parameter for each population, i.e.,
\[
r_\sigma e^{i\psi_\sigma} = \frac{1}{N_\sigma} \sum_{j=1}^{N_\sigma} e^{i\sigma_j}.
\]
Here, $r_\sigma$ describes the degree of synchronization in each population, and
ranges from $0$ to $1$. Using this, the above equations can be expressed
\[
\frac{d\sigma_i}{dt} = \omega_{\sigma i}
+  \sum_{\sigma'}\eta k_{\sigma \sigma'} r_{\sigma'} sin(\psi_{\sigma'} - \sigma_i 
- \alpha_{\sigma \sigma'}).
\]

Assuming that the $N_\sigma$ are very large, we solve for the onset of
coherent collective behavior by using a distribution function
approach. Thus, instead of discrete indicies, we imagine continua of oscillators
described by distribution functions $F_\sigma(\sigma,\omega_\sigma, t)$ such that
$F_\sigma(\sigma,\omega_\sigma, t)d\sigma d\omega_\sigma$
is the fraction of $\sigma$-oscillators whose phases and natural
frequencies lie in the infinitesimal volume element
$d\sigma d\omega_\sigma$ at time t. Note that in the $N_\sigma \rightarrow \infty$ limit,
\[
G_\sigma(\omega_\sigma)=\int_{0}^{2\pi}F_\sigma(\sigma,\omega_\sigma,t)d\sigma,
\]
and the order parameters are given by
\begin{equation}
r_\sigma e^{i\psi_\sigma} = \int \int F_\sigma e^{i\sigma} d\sigma d\omega_\sigma.
\label{orderparam}
\end{equation}
In this context, the distribution functions satisfy the equations of continuity, i.e.,
\[
\frac{\partial F_\sigma}{\partial t} + \vec\nabla \cdot
(F_\sigma \frac{d\sigma}{dt}\hat{\sigma} ) = 0,
\]
and if we write $F_\sigma(\sigma,\omega_\sigma,t)=f_\sigma(\sigma,\omega_\sigma,t)
G_\sigma(\omega_\sigma)$, we have
\begin{equation}
\frac{\partial f_\sigma}{\partial t} + \frac{\partial}{\partial \sigma}
\left\{ \left[ \omega_\sigma
+ \sum_{\sigma'}\eta  k_{\sigma \sigma'} r_{\sigma'} sin(\psi_{\sigma'} - \sigma
- \alpha_{\sigma \sigma'})
\right] \\
f_\sigma \right\} = 0.
\label{continuity}
\end{equation}

The incoherent state has $\sigma$ distributed uniformly over $[0,2\pi]$, so that
$f_\sigma=1/2\pi$ and $r_\sigma=0$. These satisfy Eq.~(\ref{continuity}) trivially.
We test the stability of this
solution by perturbing it. Note that a perturbation to $f_\sigma$ leads to a perturbation of
$r_\sigma$, and we expect these perturbations to either grow or shrink exponentially in
time, depending on the overall coupling strength $\eta$. The marginally stable case
occurs at a critical value $\eta_*$ at which coherent collective behavior
emerges.

Thus, we write $f_\sigma=1/2\pi + (\delta f_\sigma) e^{st}$ and
$r_\sigma=(\delta r_\sigma) e^{st}$, where $(\delta f_\sigma)$ and $(\delta r_\sigma)$ are small.
Inserting the first of these into the
continuity equation (Eq.~(\ref{continuity})), and keeping only first-order terms, we obtain
\[
s(\delta f_\sigma) + \omega_\sigma \frac{\partial}{\partial \sigma}(\delta f_\sigma) =
\frac{1}{2\pi} \sum_{\sigma'}\eta k_{\sigma \sigma'}(\delta r_{\sigma'}) cos(\psi_{\sigma'} - \sigma
- \alpha_{\sigma \sigma'}).
\]
The solution to this equation is
\begin{equation}
(\delta f_\sigma) = \frac{1}{4\pi}
\sum_{\sigma'} \eta k_{\sigma \sigma'}(\delta r_{\sigma'}) \left[
\frac{e^{i(\psi_{\sigma'} - \sigma -\alpha_{\sigma \sigma'})}}{s-i\omega_\sigma}
+ \frac{e^{-i(\psi_{\sigma'} - \sigma-\alpha_{\sigma \sigma'})}}{s+i\omega_\sigma}
\right].
\label{fsoln}
\end{equation}

Consistency demands that the perturbations $(\delta f_\sigma)$ and
$(\delta r_\sigma)$ be related to each other via the order parameter equation,
Eq.~(\ref{orderparam}). This yields our main result, as follows.
Eq.~(\ref{orderparam}) becomes
\[
(\delta r_\sigma)e^{st}e^{i\psi_\sigma}=
\int_{-\infty}^{\infty} G_\sigma(\omega_\sigma)
\int_{0}^{2\pi} \left( \frac{1}{2\pi} + (\delta f_\sigma)e^{st} \right)
e^{i\sigma} d\sigma d\omega_\sigma.
\]
The integral involving $1/2\pi$ is zero. Inserting the solution for $(\delta f_\sigma)$
from Eq.~(\ref{fsoln}), one obtains \cite{analyticnote}
\[
(\delta r_\sigma)e^{i\psi_\sigma}=
\left[
\frac{1}{2} \int_{-\infty}^{\infty} \frac{G_\sigma(\omega_\sigma)}{s-i\omega_\sigma}d\omega_\sigma
\right]
\sum_{\sigma'} \eta k_{\sigma \sigma'}(\delta r_{\sigma'}) e^{i(\psi_{\sigma'}-\alpha_{\sigma \sigma'})}.
\]
Define the bracketed expression to be $g_\sigma(s)$, and define
$z_\sigma = (\delta r_\sigma)e^{i\psi_\sigma}$. Then, we have
\[
z_\sigma=g_\sigma(s) \sum_{\sigma'} \eta k_{\sigma \sigma'}e^{-i\alpha_{\sigma \sigma'}}z_{\sigma'}.
\]
Now define the complex quantity $\bar{k}_{\sigma \sigma'}= \eta k_{\sigma \sigma'} e^{-i\alpha_{\sigma \sigma'}}$. Using
the Kronecker delta $\delta_{\sigma \sigma'}$, the above equation can be written
\[
\sum_{\sigma'} \left( \bar{k}_{\sigma \sigma'} - \frac{\delta_{\sigma \sigma'}}{g_\sigma(s)}
\right) z_{\sigma'} =0.
\]
In matrix notation for the case of two populations labeled $\theta$ and $\phi$, this is
\begin{equation}
\left(
\begin{array}{cc}
\bar{k}_{\theta \theta}-\frac{1}{g_\theta(s)} & \bar{k}_{\theta \phi} \\
\bar{k}_{\phi \theta} & \bar{k}_{\phi \phi}-\frac{1}{g_\phi(s)}
\end{array}
\right)
\left(
\begin{array}{c}
z_\theta \\
z_\phi
\end{array}
\right) = 0.
\label{twobytwoeqn}
\end{equation}

This equation has a non-trivial solution if the determinant of the matrix
is zero. This condition determines the growth rate $s$ in terms of
$\eta$, $\bar{k}_{\sigma \sigma'}$, and the parameters of the distributions
$G_\sigma(\omega_\sigma)$ (via $g_\sigma(s)$).

To illustrate the resulting behavior, we note that $g_\sigma(s)$ can be
evaluated in closed form for a Cauchy-Lorentz distribution
\begin{equation}
G_\sigma(\omega_\sigma)=\frac{\Delta_\sigma}{\pi} \cdot
\frac{1}{(\omega_\sigma - \Omega_\sigma)^2 + \Delta_\sigma^2},
\label{Lor}
\end{equation}
where $\Omega_\sigma$, the mean frequency of population $\sigma$,
and the half-width-at-half-maximum $\Delta_\sigma$ are both real,
and $\Delta_\sigma$ is positive. One obtains
\[
\frac{1}{g_\sigma(\omega_\sigma)} = 2 (s+\Delta_\sigma + i\Omega_\sigma).
\]

Using this, the determinant condition for the two-population case reduces to
\[
\left[ \bar{k}_{\theta \theta} -2(s+\Delta_\theta+i\Omega_\theta) \right]
\left[ \bar{k}_{\phi \phi} -2(s+\Delta_\phi+i\Omega_\phi) \right ]
-\bar{k}_{\theta \phi}\bar{k}_{\phi \theta}=0.
\]

For simplicity, we set the phase angles $\alpha_{\sigma \sigma'}$
to zero for the remainder of this paper \cite{chimera}. In this case, the matrix elements
$\bar{k}_{\sigma \sigma'}$ are purely
real, so that $\bar{k}_{\sigma \sigma'}=\eta k_{\sigma \sigma'}$.
The determinant condition then becomes
\begin{equation}
\left[ \eta k_{\theta \theta} -2(s+\Delta_\theta+i\Omega_\theta) \right]
\left[ \eta k_{\phi \phi} -2(s+\Delta_\phi+i\Omega_\phi) \right ]
- \eta^2 k_{\theta \phi} k_{\phi \theta}=0.
\label{dispeqn1}
\end{equation}

Notice that if $\eta=0$, then $s=-\Delta_\sigma - i\Omega_\sigma$, indicating that
the incoherent state is stable for zero coupling (since $-\Delta_\sigma$ is negative).
We imagine increasing (or decreasing)
$\eta$ until $s$ crosses the imaginary axis at a critical value $\eta_*$. At this point,
the incoherent state loses stability and coherent collective behavior emerges in
the ensemble. Thus, the critical value $\eta_*$ can be determined from the
determinant condition by setting $s=iv$, where $v$ is real (so that the perturbations are marginally
stable), and equating the real and imaginary parts of the left side of Eq.~(\ref{dispeqn1})
to zero. This results in two equations which can be solved
simultaneously for the two (real) unknowns $\eta$ and $v$.

For our first example, we choose two identical populations; we set $\Delta_\theta=\Delta_\phi=\Delta$ and
$\Omega_\theta=\Omega_\phi=\Omega$. (A more generic example will follow.) Denoting $D=det({\bf K})$ and
$T=tr({\bf K})$,
we separate the real and imaginary parts of Eq.~(\ref{dispeqn1})
to obtain
\begin{equation}
\left.
\begin{array}{c}
D\eta^2 - 2 \Delta T\eta + 4\Delta^2 - 4(v+\Omega)^2 = 0 \\
(v+\Omega)(4\Delta-\eta T) = 0.
\end{array}
\right.
\label{ideqns}
\end{equation}

One solution to these equations is
\[
v_*=-\Omega, \; \; \;
\eta_* = \Delta \left(
\frac{T \pm \sqrt{T^2-4D}}{D}
\right),
\]
which is valid for $T^2\geq4D$, since $\eta$ is assumed
to be real. Notice that the appropriate solution as $D \rightarrow 0$
(using the negative sign for $T>0$ and the positive sign for $T<0$)
is 
\[
v_*=-\Omega, \; \; \; \eta_*=\frac{2\Delta}{T},
\]
as can be verified by setting $D=0$ in Eq.~(\ref{ideqns}) directly.
Another solution is
\[
v_*=-\Omega \pm \frac{\Delta}{T}\sqrt{4D-T^2}, \; \; \;
\eta_*=\frac{4\Delta}{T}.
\]
Finally, setting $T=0$ in Eq.~(\ref{ideqns}) yields
$v_*=-\Omega$ and $\eta_*=\pm\frac{2\Delta}{\sqrt{-D}}$ for $D<0$, and no
solution for $D\geq0$. These results are summarized in Table \ref{identicalsolnstable}.
\begin{table}
\begin{center}
\begin{tabular}{|c|c|c|c|}
\hline
Case & Condition & $v_*$ & $\eta_*$ \\
\hline
1 & $T^2 > 4D$ & $-\Omega$ & $\Delta \left( \frac{T \pm \sqrt{T^2-4D}}{D} \right)$ \\
\hline
2 & $T^2 \leq 4D$ & $-\Omega \pm \frac{\Delta}{T}\sqrt{4D-T^2}$ & $\frac{4\Delta}{T}$ \\
\hline
3 & $D=0$ & $-\Omega$ & $\frac{2\Delta}{T}$ \\
\hline
4 & $T=0$, $D<0$ & $-\Omega$ & $\pm\frac{2\Delta}{\sqrt{-D}}$ \\
\hline
5 & $T=0$, $D\geq0$ & no solution & no solution \\
\hline
\end{tabular}
\end{center}
\caption{Solutions to Eq.~(\ref{ideqns}) for two identical populations. $D=det({\bf K})$
and $T=tr({\bf K})$, where ${\bf K}$ is the connectivity matrix (Eq.~(\ref{Kmatrix})), and
$\Delta$ is the width parameter in Eq.~(\ref{Lor}).}
\label{identicalsolnstable}
\end{table}

\begin{figure}
\begin{center}
\scalebox{0.9}{\includegraphics{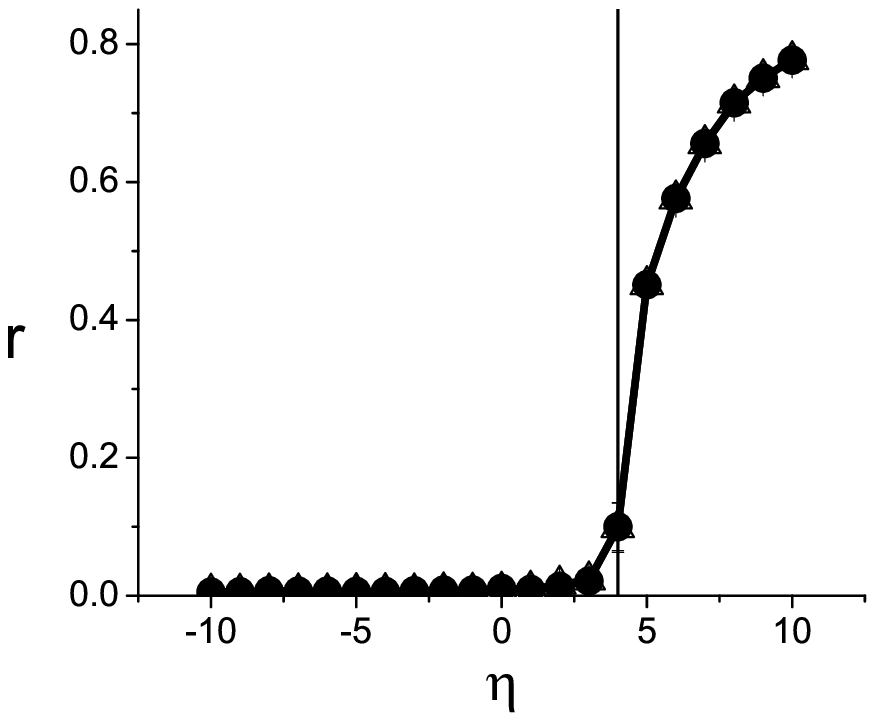}}
\caption{Numerical calculation of the order parameters ($\bullet$, $\bigtriangleup$) versus $\eta$
for case $A$ in Table~\ref{casetable}. The vertical line corresponds to the predicted value $\eta_*=4$.
The data point nearest $\eta_*$ is at $\eta=4.15$.}
\label{CaseA_samepops}
\end{center}
\end{figure}

Thus, the critical values $\eta_*$ are determined by $T$, $D$, and $\Delta$.
To illustrate this result, we begin by discussing a particular example. Consider
the matrix
\[
{\bf K}=\left(
\begin{array}{cc}
1 & -1 \\
1 & 0
\end{array}
\right),
\]
which has trace $T=1$ and determinant $D=1$. This corresponds to case
two in Table \ref{identicalsolnstable}, from which we find that $\eta_*=4\Delta/T$.
Fig.~\ref{CaseA_samepops} shows the results of a numerical simulation
of two populations of $10,000$ oscillators each, using $\Delta=1$. The order parameters $r_\sigma$ versus
$\eta$ are shown, and we can see that the oscillator populations are incoherent for $\eta$
values below the predicted critical value $\eta_*=4$, and that they grow increasingly synchronized
as $\eta$ is increased beyond $\eta_*$.

\begin{figure}
\begin{center}
\scalebox{1.7}{\includegraphics{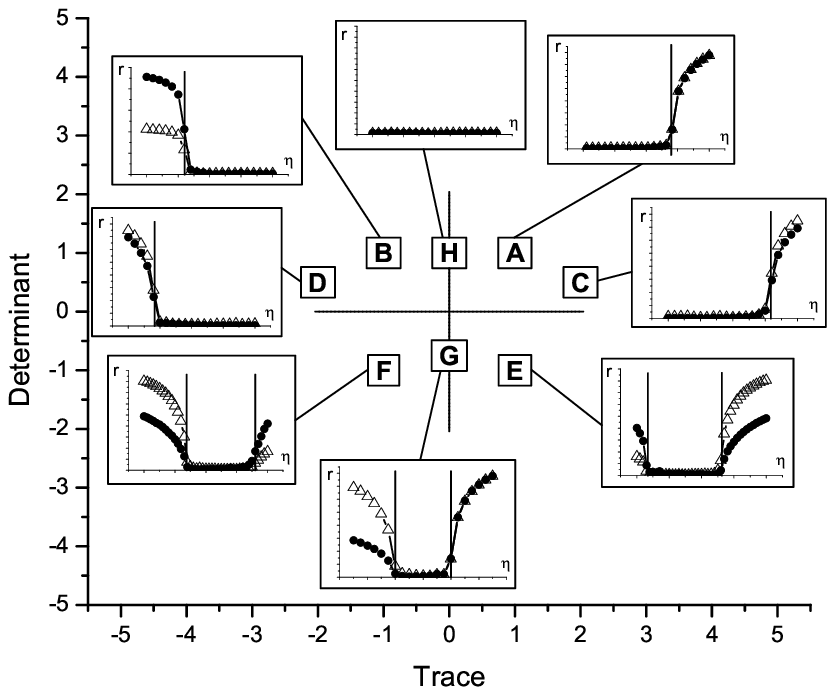}}
\caption{Numerical simulations using the matricies listed in Table (\ref{casetable})
for identical populations.
The letters indicate the placement of each case in the $T-D$ plane, and the
corresponding insets show numerical calculations of the order parameters ($\bullet$, $\bigtriangleup$) versus $\eta$
for that case (in all cases, $\eta=0$ is in the center of the horizontal axis). Vertical lines in the insets
indicate the predicted value(s) $\eta_*$ listed in Table (\ref{casetable}) for the onset of coherent collective behavior.
In all cases, we used $\Delta=1$. Note that for $D>0$, there is one value of $\eta_*$, whose sign corresponds
to the sign of the trace $T$. If $D \geq 0$ and $T=0$, then synchronization is not possible for any coupling
strength. For $D<0$, there are two values of $\eta_*$: one positive, and one negative.}
\label{masterfig}
\end{center}
\end{figure}

Next, we examined eight different connectivity matricies ${\bf K}$
that were chosen to sample
the various regions in $T-D$ space. Table~(\ref{casetable}) shows these
matricies and the corresponding value(s) of $\eta_*$.
\begin{table}
\begin{center}
\begin{tabular}{|c|c|r|r|c|}
\hline
Case & Matrix & T & D & $\eta_*$ \\
\hline
A & $\left( \begin{array}{cc} 1 & -1 \\ 1 & 0 \end{array} \right)$ & $1$ & $1$ & $4$ \\
\hline
B & $\left( \begin{array}{cc} -2 & -3 \\ 1 & 1 \end{array} \right)$ & $-1$ & $1$ & $-4$ \\
\hline
C & $\left( \begin{array}{cc} 3 & 1 \\ -3.5 & -1 \end{array} \right)$ & $2$ & $0.5$ & $2(2-\sqrt{2})=1.172$ \\
\hline
D & $\left( \begin{array}{cc} -3 & 1 \\ -3.5 & 1 \end{array} \right)$ & $-2$ & $0.5$ & $2(-2+\sqrt{2})=-1.172$ \\
\hline
E & $\left( \begin{array}{cc} -1 & -1 \\ 1 & 2 \end{array} \right)$ & $1$ & $-1$ & $-(1 \pm \sqrt{5})=-3.236, 1.236$ \\
\hline
F & $\left( \begin{array}{cc} 1 & 1 \\ -1 & -2 \end{array} \right)$ & $-1$ & $-1$ & $1 \pm \sqrt{5}=-1.236, 3.236$ \\
\hline
G & $\left( \begin{array}{cc} 2 & 1 \\ -3 & -2 \end{array} \right)$ & $0$ & $-1$ & $\pm 2$ \\
\hline
H & $\left( \begin{array}{cc} 1 & -1 \\ 2 & -1 \end{array} \right)$ & $0$ & $1$ & None \\
\hline
\end{tabular}
\end{center}
\caption{Connectivity matricies ${\bf K}$ chosen to sample $T-D$ space.}
\label{casetable}
\end{table}
These predictions were tested by numerically calculating the order
parameters $r_\sigma$ for a range of coupling values $\eta$ \cite{methods}.
Results are shown in Fig.~\ref{masterfig}. The various cases are located by letter
in the $D-T$ plane according to the trace and determinant of the matricies,
and the corresponding inset shows the numerically-calculated order parameters plotted
versus $\eta$, with the predicted critical coupling indicated by a vertical
line at $\eta=\eta_*$. For example, the inset to case $A$ (with $(T,D)=(1,1)$) corresponds to
Fig.~\ref{CaseA_samepops}, which was discussed above.
It can be seen that for all cases, the onset of synchronization occurs as predicted.

Note that more than one prediction for $\eta_*$ may be specified by our analysis
(see Table \ref{identicalsolnstable}).
The solutions closest to $\eta=0$ are the relevant ones, because we expect
the incoherent state to lose stability once the
first $\eta_*$ solution is encountered. There are two possible cases depending on
the sign of $D$. First, if the two solutions
have the same sign, then there is only one critical value $\eta_*$ (which
may be positive or negative depending on the sign of the trace) for the onset of synchronization. This occurs
for $D>0$ and $T \neq 0$, as can be seen in Fig.~\ref{masterfig}. (Interestingly,
for $D>0$ and $T=0$, synchronization does not occur for {\it any} $\eta$.) The other case occurs
for $D<0$, for which the two $\eta_*$ solutions have opposite signs.
In this case, there are two critical values $\eta_*$ for the onset of
synchronization -- one on either side of $\eta=0$. This can also be observed in
Fig.~\ref{masterfig}.

In the more general case in which the various populations have different
natural frequency distributions, it is
not typically possible to describe the onset of synchronization in terms
of the determinant and trace of the coupling matrix ${\bf K}=k_{\sigma\sigma'}$ alone.
We now consider this situation, but retain the Cauchy-Lorentz form of the
natural frequency distributions for convenience.
We manipulate Eq.~(\ref{dispeqn1}) as follows. Let $s=iv$ (i.e., purely
imaginary, to consider
the marginally stable case) and define $a=\Delta_\theta + i(v+\Omega_\theta)$
and $b=\Delta_\phi + i(v+\Omega_\phi)$. We obtain
\begin{equation}
D\eta^2 - 2\eta(bk_{\theta\theta}+ak_{\phi\phi})+4ab=0.
\label{extwoeqns}
\end{equation}
There are two unknowns in this equation: $\eta$ and $v$. Eq.~(\ref{extwoeqns}) is a quadratic
equation in $\eta$ with complex coefficients, and we can easily obtain two {\it complex} solutions $\eta^{(1,2)}$
as functions of $v$. Since the critical values $\eta_*^{(1,2)}$ must be real, we solve for the roots of $Im(\eta^{(1,2)})$,
and evaluate $Re(\eta^{(1,2)})$ at these roots. This yields the possible critical values $\eta_*^{(1,2)}$. Typically, these steps must be performed
symbolically and/or numerically;
we used MATLAB$^{\textregistered}$
\cite{matlab}.
As before, the values of $\eta_*^{(1,2)}$ that are closest to zero
(on either side) are the relevant values.

\begin{figure}
\begin{center}
\scalebox{0.19}{\includegraphics{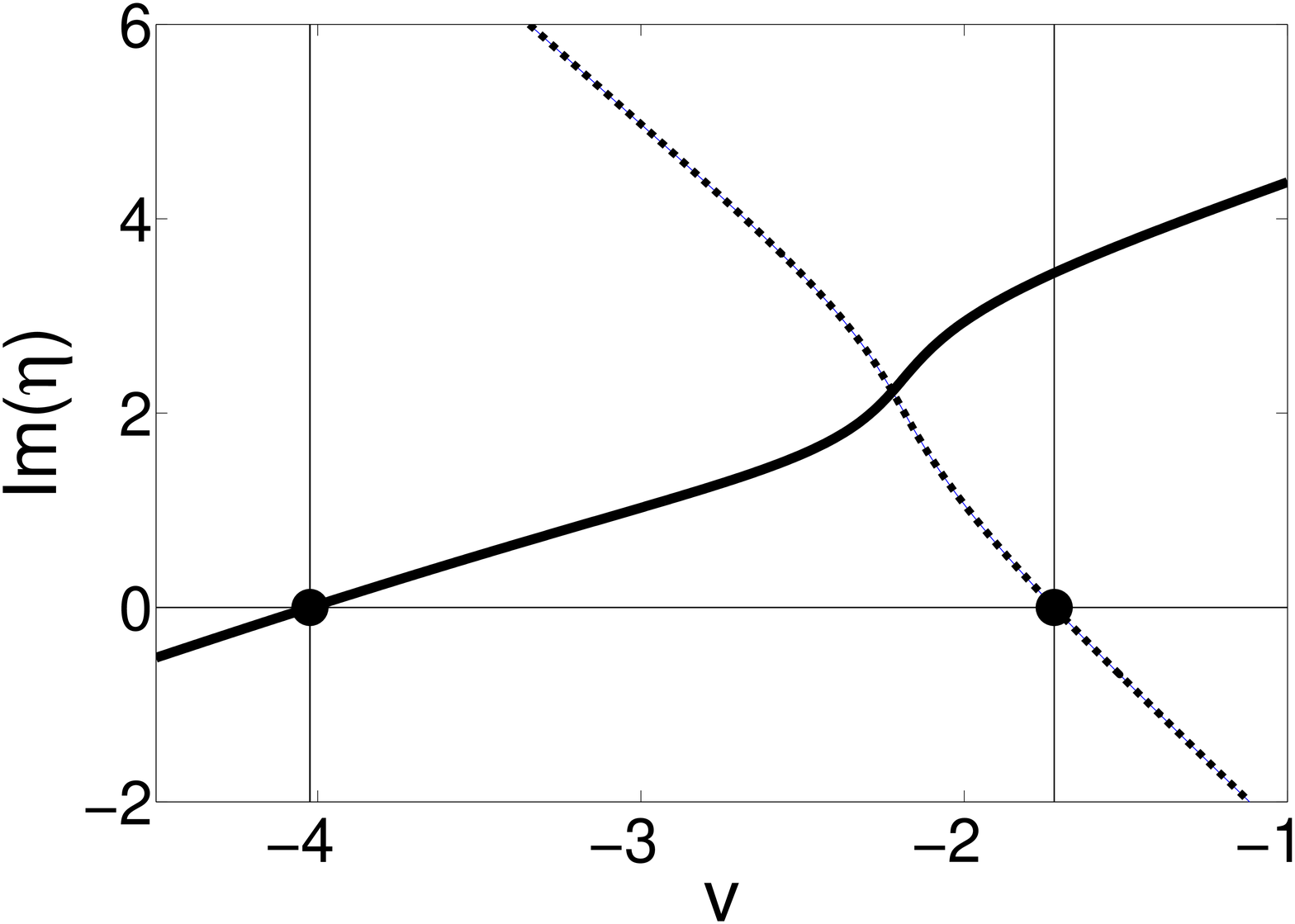}} \scalebox{0.19}{\includegraphics{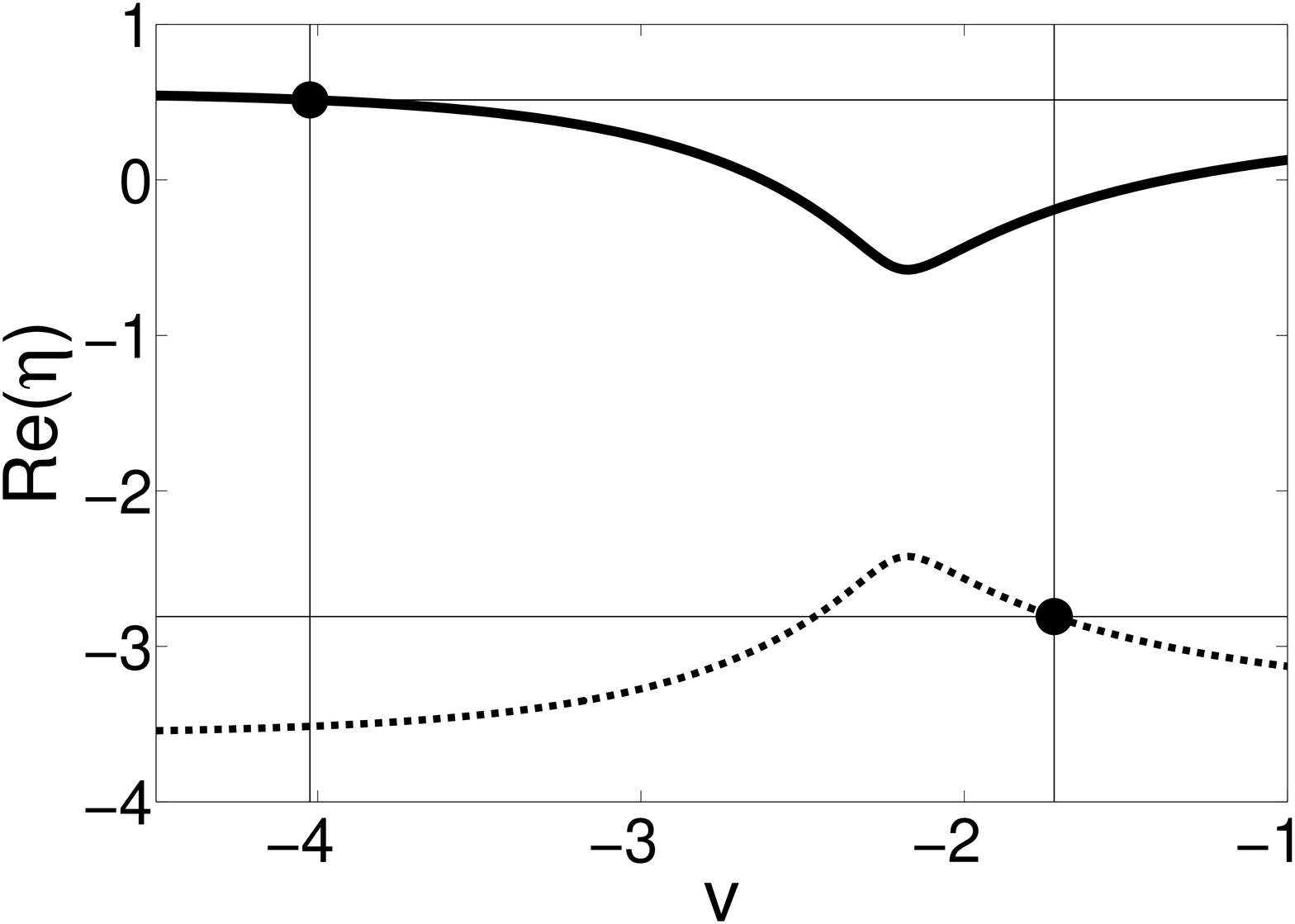}}
\caption{Case E, different populations. The left panel shows $Im(\eta^{(1,2)}$), with
the vertical lines identifying roots at $v_1=-4.024$ and $v_2=-1.722$.
The right panel shows $Re(\eta^{(1,2)})$; values at the roots found above
are indicated by horizontal lines, yielding $\eta_*^{(1)}=0.515$ and $\eta_*^{(2)}=-2.809$.
Thus, we expect
synchronization to occur at these values as $\eta$ is either increased
or decreased away from zero.}
\label{CaseEdiffpops_eta}
\end{center}
\end{figure}

\begin{figure}
\begin{center}
\scalebox{0.9}{\includegraphics{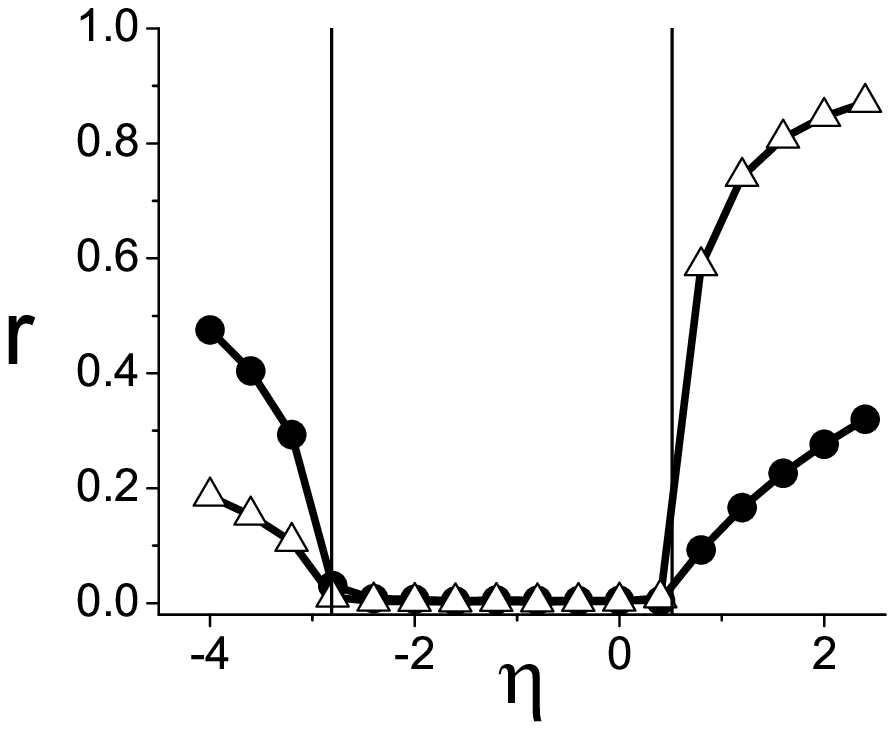}}
\caption{Case E, different populations. Calculations of the order parameters versus
$\eta$ confirm that the onset of synchronization occurs at $\eta_*=-2.809$ and $0.515$ (vertical
lines), as predicted in
Fig.~(\ref{CaseEdiffpops_eta}).}
\label{CaseEdiffpops_expt}
\end{center}
\end{figure}

\begin{figure}
\begin{center}
\scalebox{0.19}{\includegraphics{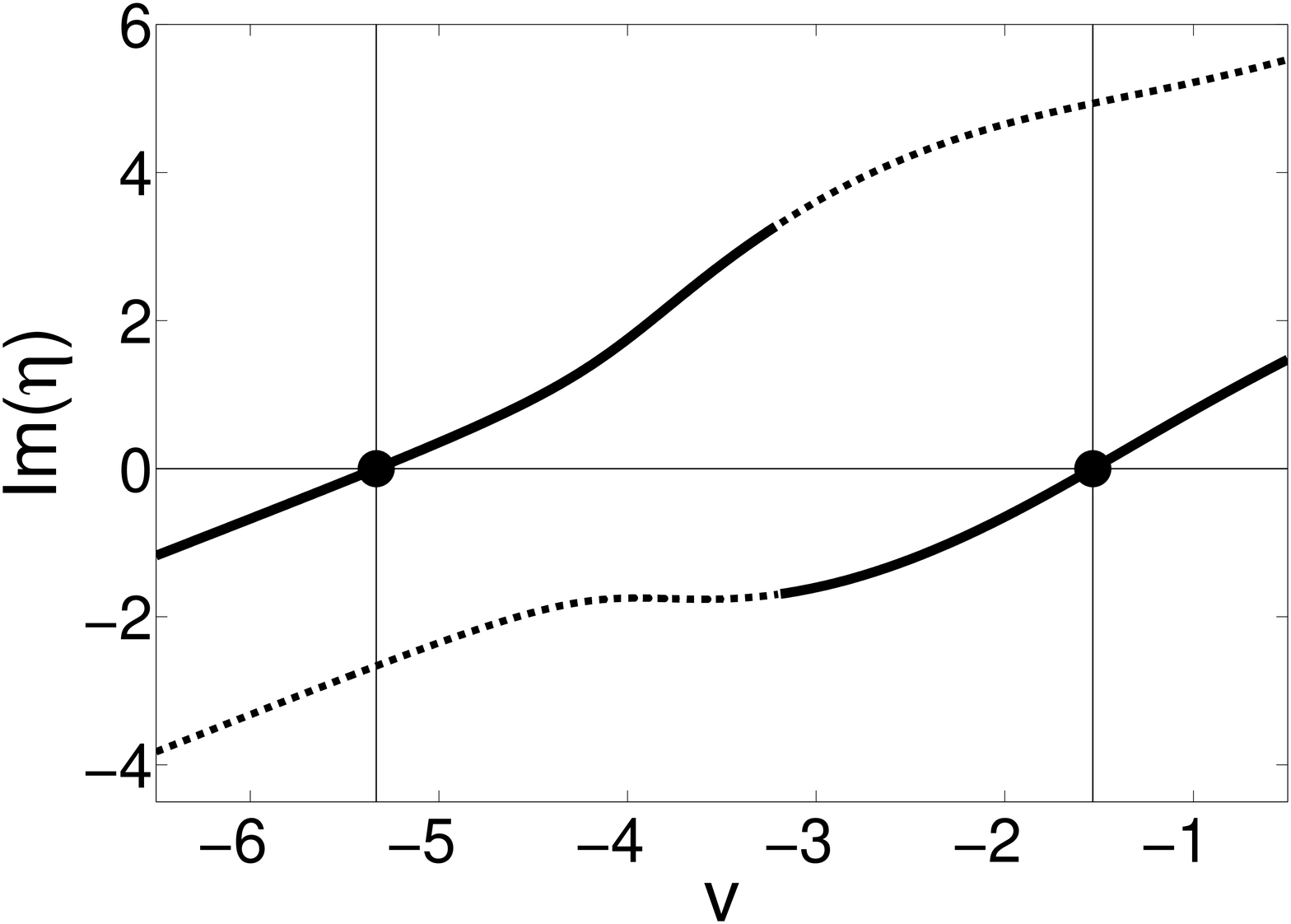}}
\scalebox{0.19}{\includegraphics{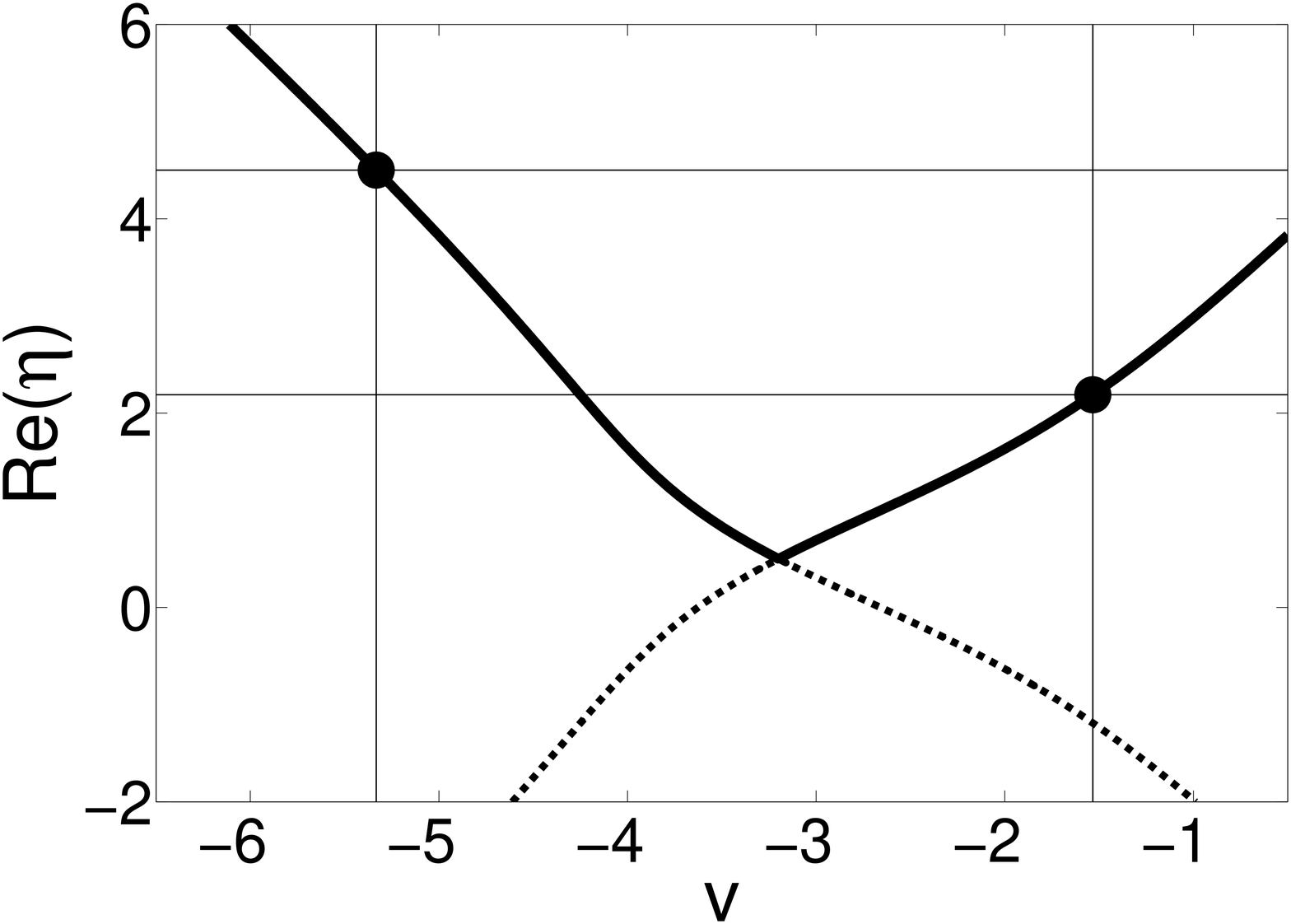}}
\caption{Case A, different populations. Panels as in Fig.~(\ref{CaseEdiffpops_eta}).
Note that the standard branch cut
leads to discontinutites and the occurrence of
two roots for $Im(\eta^{(1)})(v)$ (solid lines, left panel), and none for $Im(\eta^{(2)})(v)$
(dotted lines, left panel).
From the right panel we find $\eta_*=2.189,4.501$, taking care to obtain these from
$Re(\eta^{(1)})$ (solid line, right panel).
We expect to find synchronization onset at the smaller of these values.}
\label{CaseAdiffpops_eta}
\end{center}
\end{figure}

\begin{figure}
\begin{center}
\scalebox{0.9}{\includegraphics{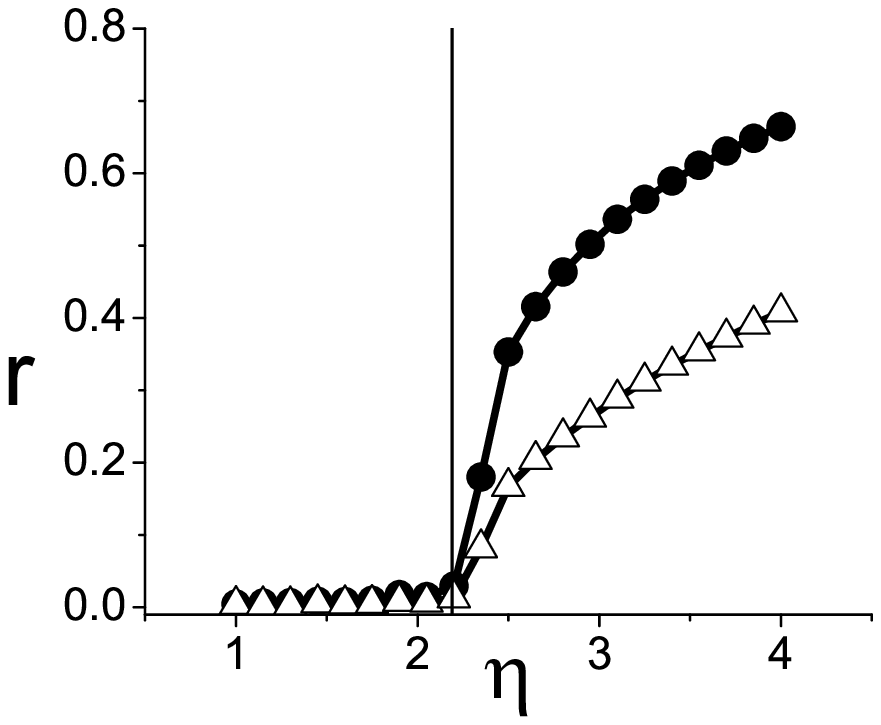}}
\caption{Case A, different populations: The onset of synchronziation occurs
at $\eta_*=2.189$ (vertical line), as predicted in Fig.~(\ref{CaseAdiffpops_eta}).
No synchronization is observed for smaller values of $\eta$ (not shown).}
\label{CaseAdiffpops_expt}
\end{center}
\end{figure}

\begin{figure}
\begin{center}
\scalebox{0.19}{\includegraphics{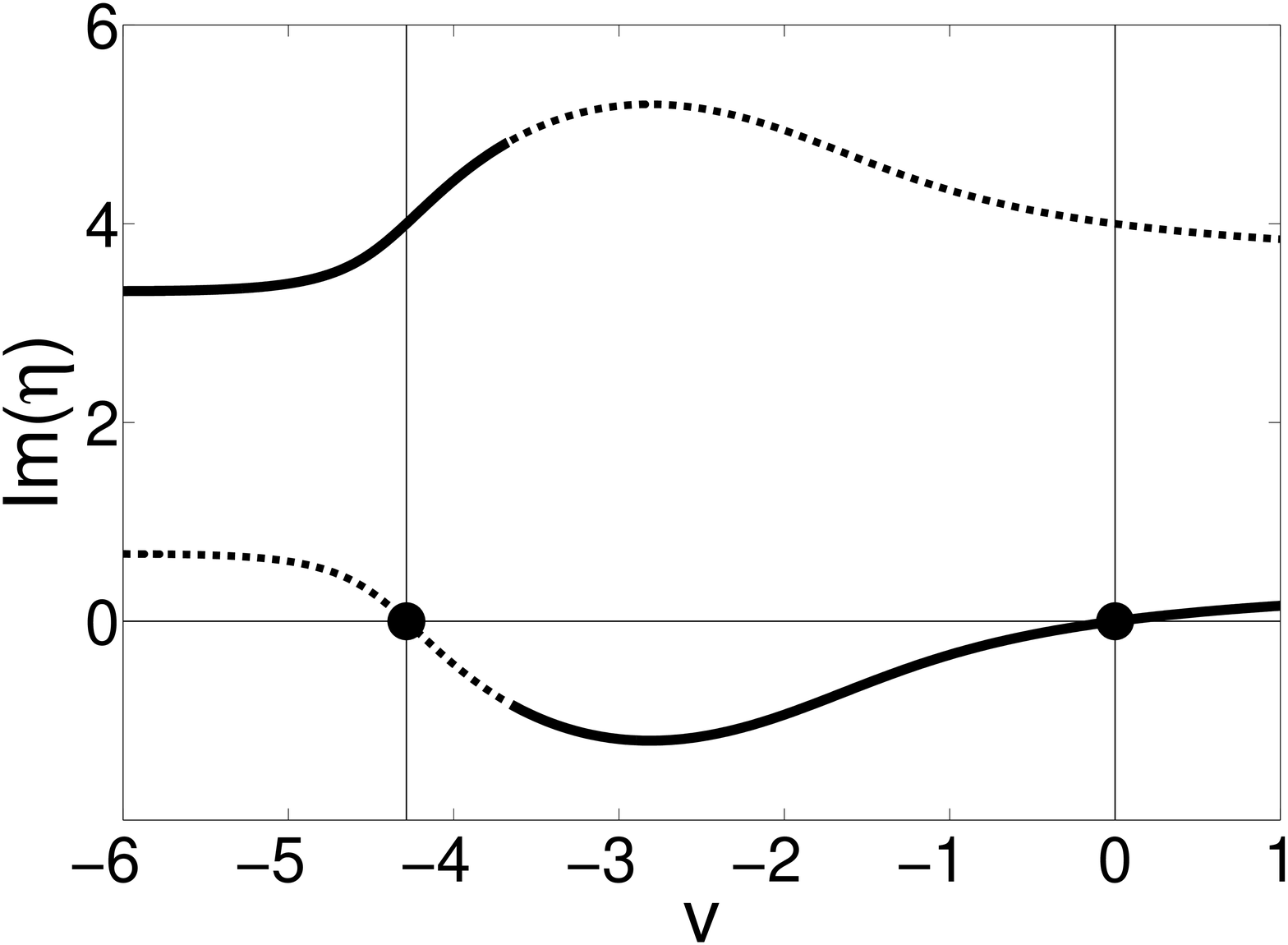}}
\scalebox{0.19}{\includegraphics{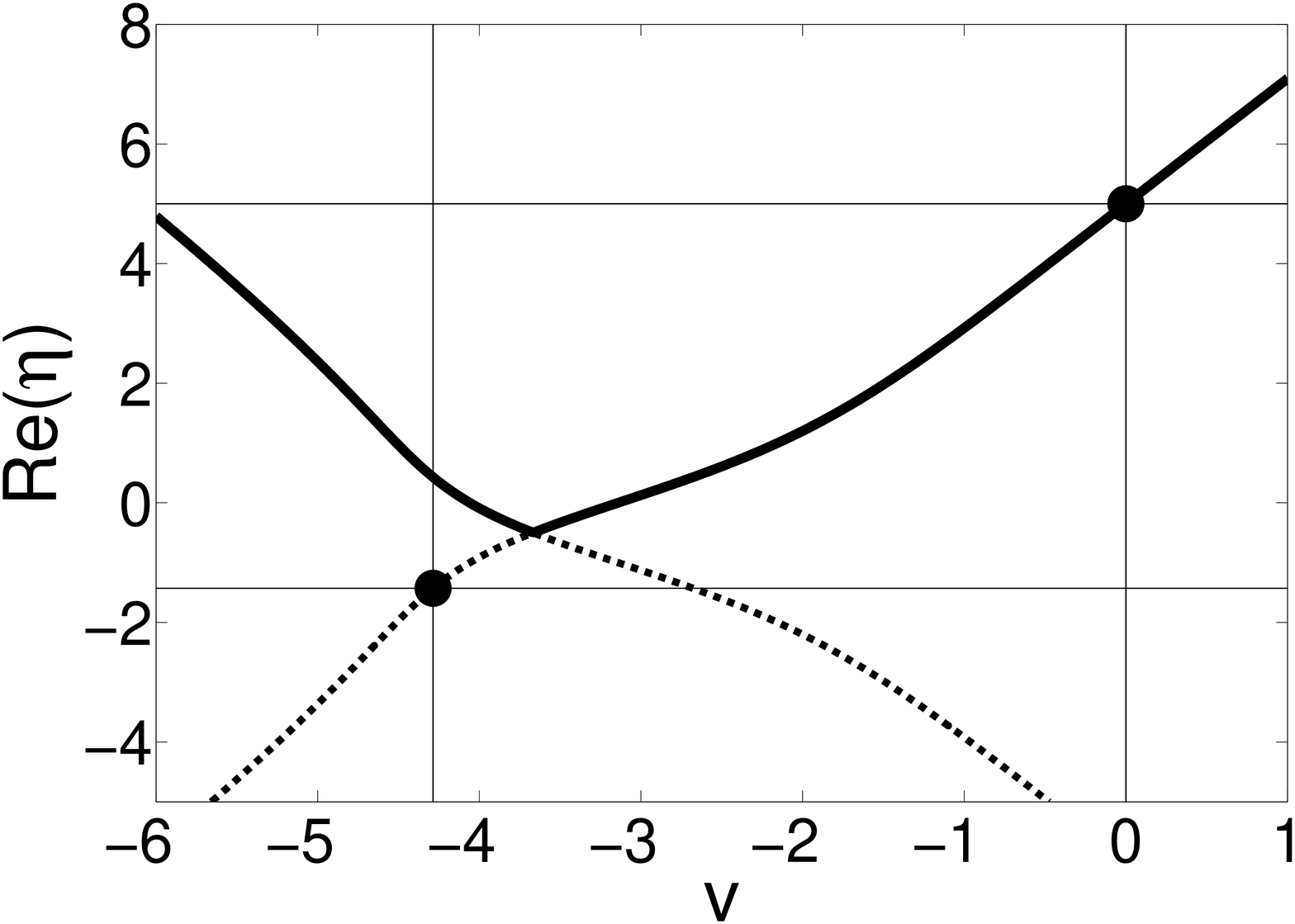}}
\caption{Case H, different populations. We find $\eta_*=-1.429,5.000$.}
\label{CaseHdiffpops_eta}
\end{center}
\end{figure}

\begin{figure}
\begin{center}
\scalebox{0.9}{\includegraphics{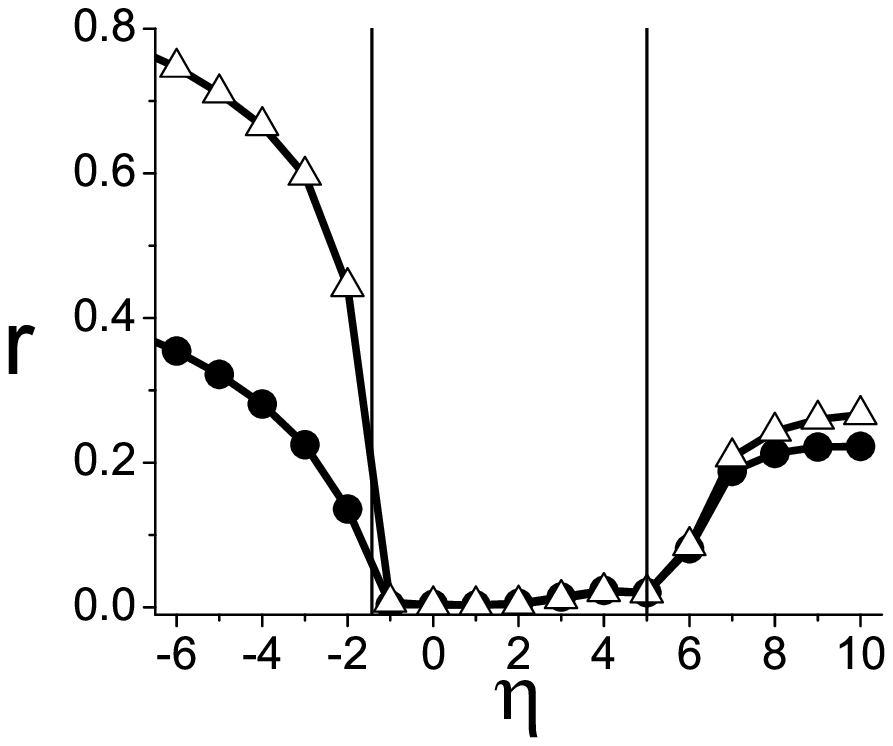}}
\caption{Case H, different populations. Synchronization occurs at $\eta_*=-1.429$
and $5.000$ (vertical lines), as predicted in Fig.~(\ref{CaseHdiffpops_eta}).}
\label{CaseHdiffpops_expt}
\end{center}
\end{figure}

To illustrate, we choose two populations with Cauchy-Lorentz natural frequency
distributions (Eq.~(\ref{Lor})) with parameters $\Delta_\theta=1$, $\Omega_\theta=2$,
$\Delta_\phi=0.5$, and $\Omega_\phi=4$.
We consider
the same $\bf{K}$ matricies as before, i.e., those listed in the second column of Table \ref{casetable}.
Case E is straightforward; the analysis is illustrated in Fig.~(\ref{CaseEdiffpops_eta})
and the numerical verification is in Fig.~(\ref{CaseEdiffpops_expt}).
Note that since Eq.~(\ref{extwoeqns}) has complex coefficients, obtaining $\eta$ typically involves
taking the square root of a complex number. Therefore, one must be mindful of branch cuts when obtaining
symbolic and/or numerical solutions. This is important in the analysis for case A, shown in
Figs.~(\ref{CaseAdiffpops_eta}) and (\ref{CaseAdiffpops_expt}).
Finally, case H, which exhibits no synchronization for identical populations
for any value of $\eta$, does show synchronization in the present case. The
analysis is shown in Figs.~(\ref{CaseHdiffpops_eta}-\ref{CaseHdiffpops_expt}).

\begin{figure}
\begin{center}
\scalebox{0.19}{\includegraphics{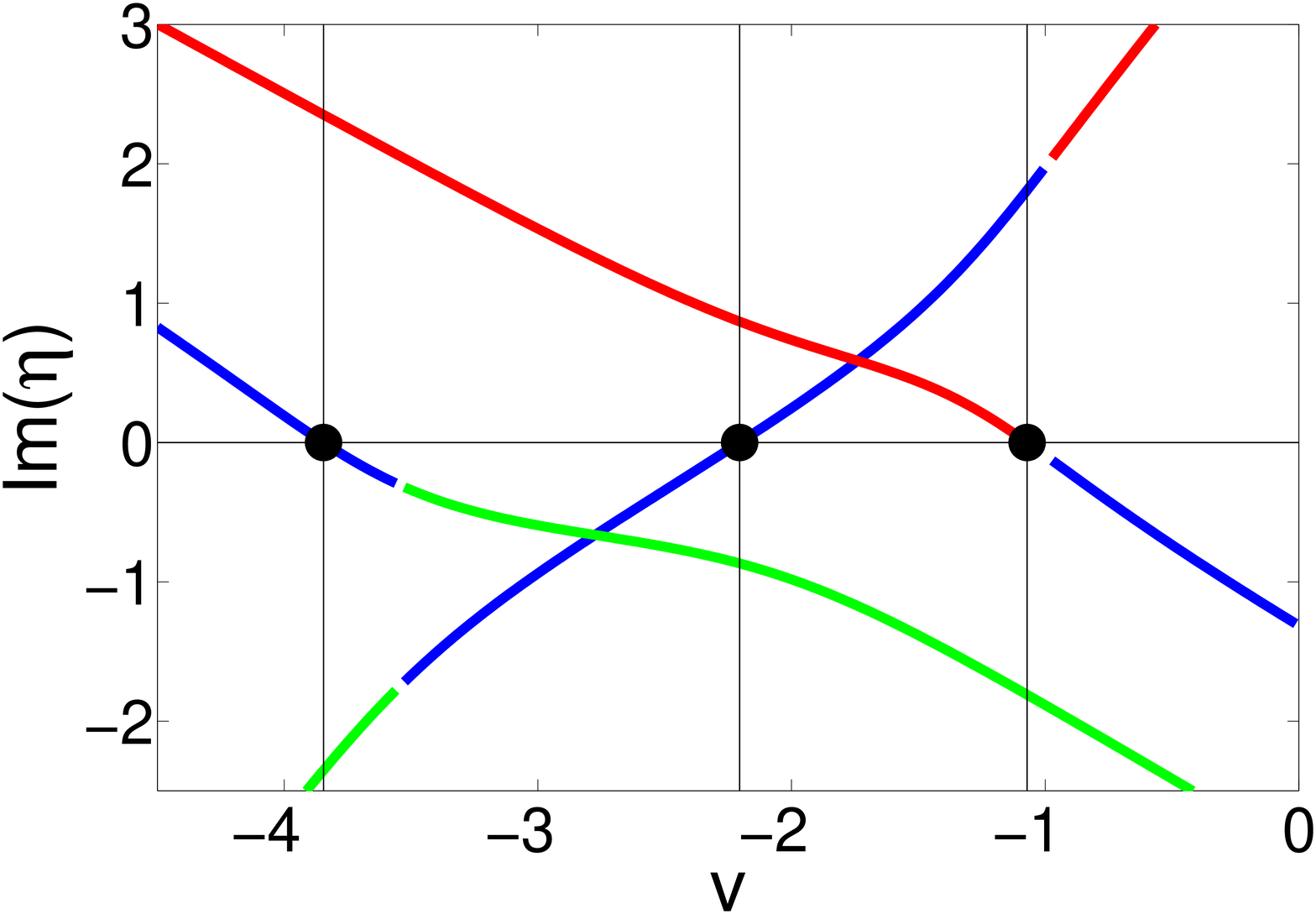}}
\scalebox{0.19}{\includegraphics{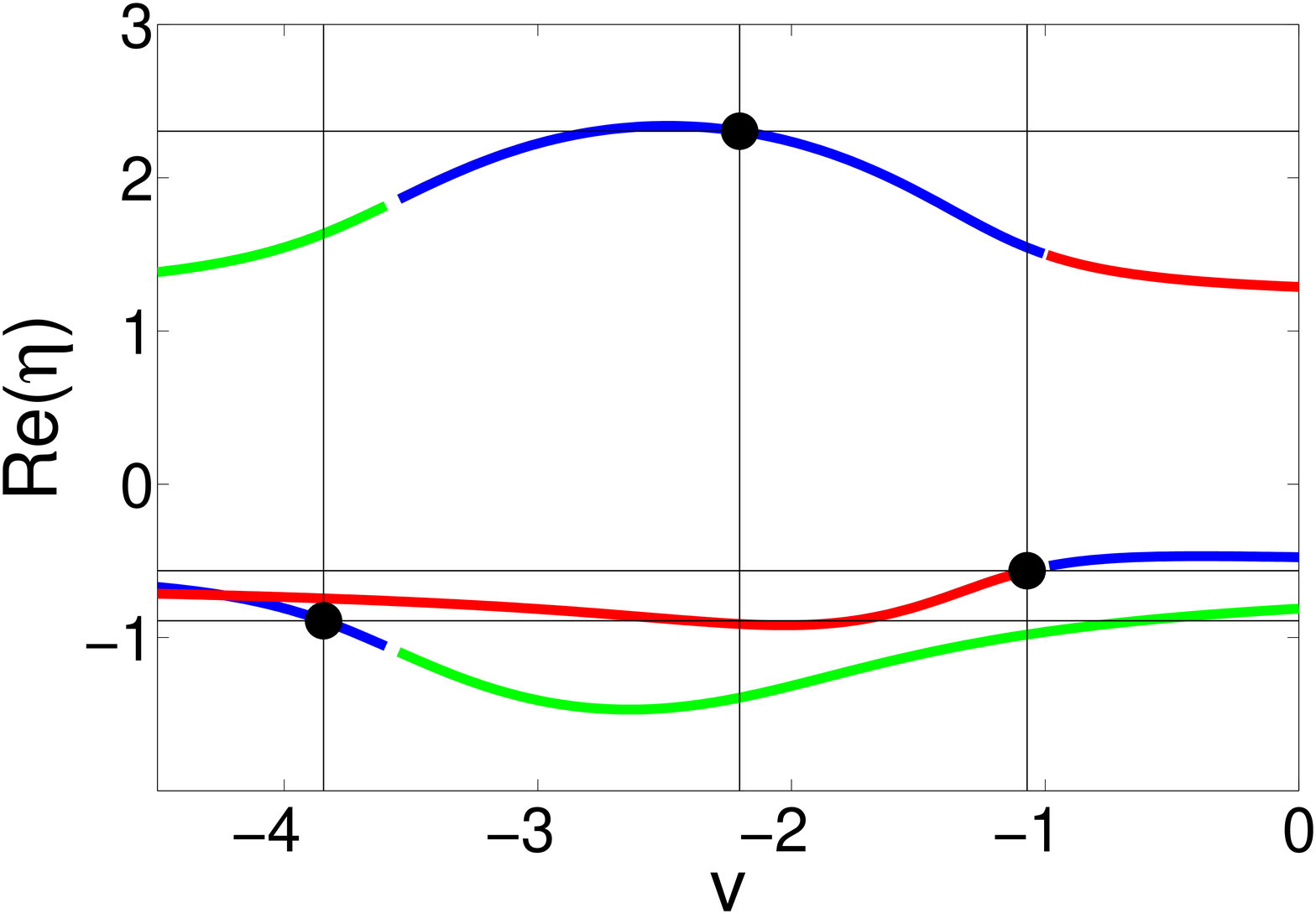}}
\caption{(Color online) Three populations. The imaginary and real parts of $\eta^{(1,2,3)}$
are plotted in blue, red, and green to illustrate the discontinuties due to branch cuts. The analysis proceeds as in
the previous cases. Because of branch cuts, two roots occur on $\eta^{(1)}$, one on $\eta^{(2)}$,
and none on $\eta^{(3)}$. Evaluating the real parts appropriately, we find $\eta_*=-0.891,-0.564,2.303$.}
\label{threepredict}
\end{center}
\end{figure}

\begin{figure}
\begin{center}
\scalebox{0.9}{\includegraphics{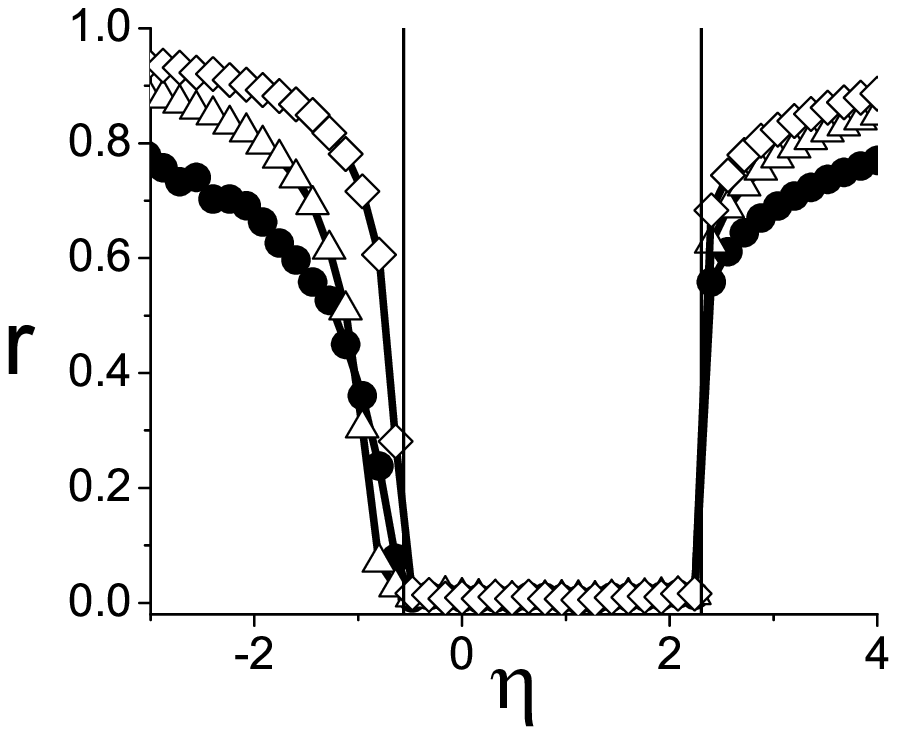}}
\caption{Three populations. Synchronization occurs at $\eta_*=-0.564$ and $2.303$,
as predicted in Fig.~(\ref{threepredict}).}
\label{threepop}
\end{center}
\end{figure}

We close by giving an example with three different populations. We choose
the same Cauchy-Lorentz distributions as above, and add a third with 
$\Delta_\rho=1/3$ and $\Omega_\rho=1$. We use the following $\bf{K}$ matrix:
\[
\left(
\begin{array}{ccc}
-1 & 1 & 1 \\
1 & -1 & 1 \\
1 & 1 & -1
\end{array}
\right).
\]
The procedure for deriving $\eta_*$ proceeds as above, except that Eq.~(\ref{dispeqn1})
is replaced by a third-degree polynomial
in $\eta$. Fig.~(\ref{threepredict}) shows the imaginary and real parts of the
three $\eta$ solutions. (Note that the branch cuts are more complicated.)
The predicted onset of synchronization was verified,
as shown in Fig.~(\ref{threepop}).

In conclusion, we have described how to determine the onset of coherent collective
behavior in systems of interacting Kuramoto systems, i.e., systems of interacting
{\it populations} of phase oscillators with both node and coupling heterogeneity.
EB was supported by NIH grant R01-MH79502; EO was
supported by ONR (Physics) and by NSF grant PHY045624.


\begin{thebibliography}{}

\bibitem{motifs}
R. Milo et al., Science {\bf 298}, 824-827 (2002).

\bibitem{community}
M.E.J. Newman and M. Girvan, Phys. Rev. E {\bf 69}, 026113 (2004).

\bibitem{community2}
A. Arenas, A. D\'{i}az-Guilera, and C.J. P\'{e}rez-Vicente,
Phys. Rev. Lett. {\bf 96}, 114102 (2006);
A. Arenas, A. D\'{i}az-Guilera, and C.J. P\'{e}rez-Vicente,
Physica D {\bf 224}, 27-34 (2006).

\bibitem{layers}
M. Kurant and P. Thiran, Phys. Rev. Lett. {\bf 96}, 138701 (2006).

\bibitem{cluster}
X. Wang, L. Huang, Y-C Lai, and C.H. Lai, Phys. Rev. E {\bf 76}, 056113 (2007).

\bibitem{ZhouKurths}
C. Zhou, L. Zemanov\'a, G. Zamora, C.C. Hilgetag, and J. Kurths, Phys. Rev. Lett. {\bf 97}, 238103 (2006).

\bibitem{GirvanNewman}
M. Girvan and M.E.J. Newman, Proc. Natl. Acad. Sci. {\bf 99} 7821-7826 (2002).

\bibitem{milo2}
R. Milo et al., Science {\bf 303}, 1538-1542 (2004).

\bibitem{kur84}
Y. Kuramoto, in {\em International Symposium on Mathematical Problems in Theoretical Physics}, edited by
H. Araki, Lecture Notes in Physics, Vol.~39 (Springer, Berlin, 1975); {\em Chemical Oscillators, Waves and Turbulence}
(Springer, Berlin, 1984).

\bibitem{apps}
A. T. Winfree, {\em The Geometry of Biological Time} (Springer, New York, 1980);
S.H. Strogatz, Physica D {\bf 143}, 1 (2000);
J.A. Acebr\'{o}n et al., Rev. Mod. Phys. {\bf 77}, 137-185.
K. Wiesenfeld and J.W. Swift, Phys. Rev. E {\bf 51}, 1020-1025 (1995).
I.Z. Kiss, Y. Zhai, and J.L. Hudson, Science {\bf 296} 1676 (2005).

\bibitem{Montbrio}
A similar system of two asymmetrically interacting populations with a particular form of coupling
matrix ${\bf K}$ was considered in E. Montbrio, J. Kurths, and B. Blasius, Phys. Rev. E 70, 056125 (2004).
The formulation in the current paper is more general in two important aspects: ${\bf K}$ is
arbitrary, and we allow for any number of interacting populations.

\bibitem{restreponote}
Our system is similar to that studied in J.G. Restrepo, E. Ott, and J.G. Restrepo,
Chaos {\bf 16}, 015107 (2005). However, our
formulation permits different natural frequency distributions for each population.

\bibitem{analyticnote}
The bracketed expression is valid for $Re(s)>0$ and may be analytically
continued into the region where $Re(s) \leq 0$. See E. Ott, P. So, E. Barreto,
and T. Antonsen, Physica D {\bf 173}, 226-258 (2002); E. Ott, {\it Chaos in Dynamical
Systems}, second edition, Cambridge University Press, 2002, p. 240.

\bibitem{chimera}
Many interesting states require $\alpha_{\sigma \sigma'} \neq 0$, such as the
chimera state observed in
D.M. Abrams and S.H. Strogatz, Int. J. Bif. Chaos {\bf 16} \#1, 21-37 (2006).

\bibitem{methods}
Simulations used fourth-order Runge-Kutta with a timestep of $0.01$ seconds,
$N=10,000$ or $50,000$, and parameters as noted. The system was initialized in the
incoherent state and an initial transient was discarded.
The order parameters were then averaged over the subsequent $10$ seconds.
Because the standard deviations over this interval were small, no error bars
were plotted.

\bibitem{matlab}
MATLAB$^{\textregistered}$, The MathWorks, Inc., Natick, MA ({\tt http://www.mathworks.com/}).


\end{thebibliography}
\end{document}